\documentclass [11pt,twoside]   {article}
\usepackage[margin=1.5in]{geometry}
\usepackage{amsmath}
\usepackage{amsfonts}
\usepackage{amssymb}

\usepackage{hyperref}

 \newtheorem{lemma}{Lemma}[section]
  \newtheorem{theorem}[lemma]{Theorem}
\newtheorem {example}[lemma]{Example}

  \newtheorem{remark}[lemma]{Remark}

\numberwithin{equation}{section}
  \newcommand {\pf}  {\mbox{\sc\bf Proof. \,\,}}
    \newcommand {\pff}  {\mbox{\sc\bf Proof of }}

  \newcommand {\qed} {\null \hfill \rule{2mm}{2mm}}

\begin {document}
\title{{\Large{\bf Injectivity of the specialization homomorphism of elliptic curves}}}

\author
{  {Ivica Gusi\'c and Petra Tadi\'c\thanks{ The second author was supported by the Austrian Science Fund (FWF): P 24574-N26 and both authors were supported by the Croatian Science Foundation under the project no. 6422.}}  \vspace{1ex}\\
}

\date{}
\maketitle

\begin{abstract}
\noindent  Let $E:y^2=x^3+Ax^2+Bx+C$  be a nonconstant
elliptic curve over $\mathbb{Q}(t)$ with at least one nontrivial $\mathbb{Q}(t)$-rational $2$-torsion point. We describe a method for
finding  $t_0\in\mathbb Q$ for which the corresponding specialization homomorphism $t\mapsto t_0\in\mathbb{Q}$ is injective. The method can be directly extended  to elliptic curves over $K(t)$ for a number field $K$ of class number $1$, and in principal for arbitrary number field $K$.  One can  use this method to calculate the rank of  elliptic curves over $\mathbb Q(t)$ of the form as above, and to prove that given points are free generators.
In this paper we  illustrate it on some elliptic curves over $\mathbb Q(t)$ from an article by Mestre.\end{abstract}

\footnotetext{ {\it 2000 Mathematics Subject Classification.}
11G05, 14H52.\\
{\it Key words and phrases.} elliptic curve, specialization
homomorphism, number field, class number, quadratic field, cubic field, rank, Pari, Magma}
\section{Introduction}
 \label{section1}
  Let
\begin{equation}\label{jedn}
E=E(t):y^2=x^3+Ax^2+Bx+C
\end{equation}
be a nonconstant (non-isotrivial) elliptic curve over $\mathbb{Q}(t)$ , i.e.,  $E$ is not
isomorphic over $\mathbb{Q}(t)$
to an elliptic curve over $\mathbb Q$. For the sake  of simplicity we will assume that $A,B,C\in \mathbb{Z}[t]$.  It is known that the set $E(\mathbb Q(t))$ of $\mathbb Q(t)$-rational points of $E$ is finitely generated. Let $D$ denote the discriminant of the polynomial
$f(x):=x^3+Ax^2+Bx+C$. We  note that $D\in \mathbb{Z}[t]$.
Let $t_0\in\mathbb{Q}$ be such that
$D(t_0)\neq 0.$ Then by specializing $t$ to $t_0$ the specialization
$E(t_0)$ of $E(t)$ is an elliptic curve over $\mathbb{Q}$ and we have a specialization homomorphism
 $\sigma=\sigma_{t_0}: E(\mathbb{Q}(t))\rightarrow
E(t_0)(\mathbb{Q})$ (note that it is well defined). For more on this topic see \cite[Appendix C \S 20]{Sil234}. The specialization homomorphism can be defined for general
non-split elliptic surfaces and in a more general situation. In 1952 A. N\' eron \cite{Neron} showed that the specialization fails to be injective for $t_0\in \mathbb Q$ on a small subset (of density 0) (see \cite[Section 11.1]{Serre}).  J. H. Silverman \cite{Sil0,Sil1} in 1983 using heights and J. Top  in 1985 in his master's thesis (see \cite{Top})
by extending N\'eron's techniques proved  the  so called Silverman specialization theorem, which says that the specialization homomorphism
is in fact injective for all but finitely many
rational $t_0$.  As far as
we know, there is no practical algorithm for determining such a $t_0$ (for general non-split elliptic surfaces).
As we learned from J. H. Silverman, all constants  in \cite{Sil1}, Section $4$, Theorem B, and Section $5$, Theorem $C$ can, in principal, be effectively computed. Therefore, one can find a computable constant $C$, such that for all algebraic $t_0$ with height greater than $C$, the specialization homomorphism at $t_0$ is injective. However, the constants are too large to be practical. Similarly  for methods from \cite{Sil2}.
In this
paper we use the ideas  from N\'eron and Top (which also appear in \cite{hazama}). We obtain a method for
finding a specialization $t\mapsto t_0\in\mathbb{Q}$ such that the
specialization homomorphism is injective, in the case of elliptic
curves of shape \eqref{jedn} having at least one non-trivial $\mathbb{Q}(t)$-rational $2$-torsion point. This improves and extends the method from \cite{G-T}.
Let us state the main results (see Section \ref{section2} and Section \ref{section3} for the proofs):

\begin{theorem}\label{main}
Let  $E$ be a nonconstant elliptic curve over $\mathbb Q(t)$, given by the equation
$$E=E(t):y^2=(x-e_1)(x-e_2)(x-e_3), (e_1,e_2,e_3\in \mathbb Z[t]).$$
Assume that $t_0\in\mathbb Q$ satisfies the following condition.\\
(A) For every nonconstant square-free divisor $h$ in  $\mathbb Z[t]$ of
$$\mbox{$(e_1-e_2)\cdot (e_1-e_3)$\ \  or \ \ $(e_2-e_1)\cdot(e_2-e_3)$ \ \ or \ \ $(e_3-e_1)\cdot(e_3-e_2)$},$$
the rational number $h(t_0)$ is not a square in $\mathbb{Q}$.\\
Then the specialization homomorphism $\sigma:E(\mathbb{Q}(t))\rightarrow
E(t_0)(\mathbb{Q})$ is injective.
\end{theorem}

This leads to a
 practical criterion that can be directly extended to  number fields $K$ of class number one, where the elliptic curves are as in Theorem \ref{main}  with $e_j\in \mathcal{O}_K[t]$ (here $\mathcal{O}_K$ is the ring of integers of $K$). The criterion can be extended to arbitrary number fields. However, the calculations over general number fields are rather complicated. For example, if the class number of the field $K$ is greater then $1$, the ring of integers $\mathcal{O}_K$ has to be replaced by a suitable UFD.

\begin{theorem}\label{mainh}
Let $K$ be a number field. Let $\mathcal R_K$ be  a chosen  unique factorization domain such that $\mathcal O_K\subset \mathcal R_K\subset K$ (and such that its group of units is finitely generated). If $K$ is of class number one we always choose $\mathcal R_K=\mathcal O_K$.  Let  $E$ be a nonconstant elliptic curve over $K(t)$, given by the equation
\begin{equation}\label{geneq}
E=E(t):y^2=(x-e_1)(x-e_2)(x-e_3), (e_1,e_2,e_3\in \mathcal R_K[t]).
\end{equation}
Assume that $t_0\in K$ satisfies the following condition.\\
(C) For every nonconstant square-free divisor $h$ in  ${\mathcal {R}}_K[t]$ of
$$\mbox{$(e_1-e_2)\cdot (e_1-e_3)$\ \  or \ \ $(e_2-e_1)\cdot(e_2-e_3)$ \ \ or \ \ $(e_3-e_1)\cdot(e_3-e_2)$},$$
the algebraic number $h(t_0)$ is not a square in $K$.\\
Then the specialization homomorphism $\sigma:E(K(t))\rightarrow
E(t_0)(K)$ is injective.
\end{theorem}

The first author extends Theorem \ref{main} to elliptic curves of shape \eqref{jedn} having exactly one non-trivial $\mathbb{Q}(t)$-rational $2$-torsion point. Let us set $f(x):=x^3+Ax^2+Bx+C=(x-e_1)(x-e)(x-\bar e)$ where $e_1\in \mathbb Z[t]$ and $e,\bar e$ are conjugate. Here $D=(e_1^2-(e+\bar e)e_1+e\bar e)^2(e-\bar e)^2$.

\begin{theorem}\label{novi}
Let $E=E(t):y^2=x^3+Ax^2+Bx+C;\ A,B,C\in \mathbb{Z}[t]$  be an nonconstant elliptic curve over $\mathbb{Q}(t)$.
Assume that $E$ has
 exactly one nontrivial $2$-torsion point over $\mathbb{Q}(t)$, i.e. that $x^3+Ax^2+Bx+C=(x-e_1)(x-e)(x-\bar e)$ where $e_1\in \mathbb Z[t]$ and $e,\bar e$ are conjugate. Let $t_0\in \mathbb{Q}$  satisfy the following condition:\\
 ($\mathcal{A}$) For every nonconstant square-free divisor $h$ of  $(e_1^2-(e+\bar e)e_1+e\bar e)$ or $(e-\bar e)^2$ in $\mathbb Z[t]$
the rational number $h(t_0)$ is not a square in $\mathbb{Q}$.\\
Then the specialized curve $E_{t_0}$ is elliptic and the specialization homomorphism at $t_0$ is injective.
\end{theorem}

This theorem also can be directly extended to  number fields $K$ of class number one. It can  also be extended to arbitrary number fields by passing from the ring of integers  $\mathcal{O}_K$ to a suitable ring $\mathcal R_K$ as above.

\begin{theorem}\label{mainnovi}
Let $K$ be a number field. Let $\mathcal R_K$ be  a chosen  unique factorization domain such that $\mathcal O_K\subset \mathcal R_K\subset K$ (and such that its group of units is finitely generated). If $K$ is of class number one we always choose $\mathcal R_K=\mathcal O_K$.  Let  $E$ be a nonconstant elliptic curve over $K(t)$, given by the equation
 $E:y^2=x^3+Ax^2+Bx+C;\ A,B,C\in \mathcal R_K[t]$.
Assume that $E$ has
 exactly one nontrivial $2$-torsion point over $K(t)$, i.e. that $x^3+Ax^2+Bx+C=(x-e_1)(x-e)(x-\bar e)$ where $e_1\in \mathcal R_K[t]$ and $e,\bar e$ are conjugate. Let $t_0\in K$  satisfy the following condition:\\
 ($\mathcal{A}$) For every nonconstant square-free divisor $h$ of  $(e_1^2-(e+\bar e)e_1+e\bar e)$ or $(e-\bar e)^2$ in $R_K[t]$
the  number $h(t_0)$ is not a square in $K$.\\
Then the specialized curve $E_{t_0}$ is elliptic and the specialization homomorphism at $t_0$ is injective.
\end{theorem}

The criterion from these theorems has no direct extension to general
elliptic curves  (see Remark \ref{nijeopci}).

 There are a few methods for calculating the rank of specific
types of elliptic curves over $\mathbb Q(t)$ such as the method
based on the Tate-Shioda formula (see \cite{Sh1}, Corollary 15,
and \cite{Sh2}), Nagao's conjectural method \cite{Na} and the
2-descent method (see, for example \cite{Br}).  Our method helps
in calculating  the rank of elliptic curves over $\mathbb Q(t)$
by using information about a suitably chosen specialized curve. One
can use it even to prove that a set of points are free generators.
In Section \ref{section4} we describe and comment on a family of
quadratic twists coming from Mestre: a family of quadratic twists
of the general family of elliptic curves $E=E^{a,b}:y^2=x^3+ax+b$
over $\mathbb{Q}$ with certain $14$th degree polynomials
$g=g^{a,b}$ in a variable $t$ over $\mathbb{Q}$.
 It is known that the rank of $E_g$ over $\mathbb{Q}(t)$ is at least $2$ for all $a,b,\ ab\neq 0$. By a general principle, these ranks are at most $6$. We performed an extensive calculation using our criterion  for  number fields of class number one (including $\mathbb Q$) and   for a number field of class number two. In all cases we get that the rank is two and prove that the given points are free generators. Two examples are presented in details, one when the spitting field of the polynomial $x^3+ax+b$ is a sextic field with class number one (see Example \ref{corall}), the other  when the splitting field has class number two (see Example \ref {corh2}).
We used Magma \cite{MAGMA}, Pari \cite{Pari}, and {\tt mwrank} \cite{mwrank} for most of our computations.

We would like to thank Andrej Dujella and Joseph H. Silverman for their kind suggestions and comments. We especially would like to thank the referees for very useful comments which enabled a significant improvement of the first version of the manuscript.

\section{Elliptic curve with rational $2$-torsion}
\label{section2}
In this section we  prove Theorem \ref{main} and sketch a proof of Theorem \ref{mainh}. First we work over $\mathbb{Q}$ as the field of constants. At the end of the section we extend the consideration to arbitrary number fields.
Let  $E$ be a nonconstant elliptic curve over $\mathbb Q(t)$, given by the equation
$$E=E(t):y^2=(x-e_1)(x-e_2)(x-e_3), (e_1,e_2,e_3\in \mathbb Z[t]).$$ We have homomorphisms
$\Theta_i:E(\mathbb{Q}(t))\rightarrow \mathbb{Q}(t)^{\times}/(\mathbb{Q}(t)^{\times})^2,\ i=1,2,3$  given by\\
\[\left\{\begin{array}{ll}
\Theta_i(x,y)=(x-e_i)\cdot (\mathbb{Q}(t)^{\times})^2,&\mbox{ if $x\neq e_i$,}\\
\Theta_i(e_i,0)=(e_j-e_i)(e_k-e_i)\cdot (\mathbb{Q}(t)^{\times})^2,&\mbox{ where $i\neq j\neq k\neq i$,}\\
\Theta_i(O)=1\cdot (\mathbb{Q}(t)^{\times})^2,&\mbox{ (here $O$ denotes the neutral element)}.
\end{array}
\right.\]
\begin{lemma}\label{dupla}
  $P\in 2E(\mathbb{Q}(t))$ if and only if $\Theta_i(P)=1\cdot (\mathbb{Q}(t)^{\times})^2$ for $i=1,2,3$.
\end{lemma}
\pf The claim follows from  \cite{Huse},
Chapter 1, Theorem (4.1), and  Chapter 6,
Proposition (4.3). See also \cite[Section 1]{hazama} for the generalization to hyperelliptic curves.\qed\medskip

Since $\mathbb{Z}[t]$ is a unique factorization domain (UFD), it is evident that for each $P\in E(\mathbb{Q}(t))$ there exists exactly one
triple $(s_1,s_2,s_3),\  s_i=s_i(P)\in \mathbb{Z}[t], \ i=1,2,3 $, of
non-zero square-free elements from $\mathbb{Z}[t]$, such that
\begin{equation}\label{esi}
\Theta_i(P)= s_i(P)\cdot (\mathbb{Q}(t)^{\times})^2.
\end{equation}
 We will also use notation  $s_i(t)$ for $s_i$. Lemma \ref{dupla} can be reformulated as

\begin{equation}\label{opetdupla}
P\in 2E(\mathbb{Q}(t))\ {\rm if\ and\ only\ if}\ s_i(P)=1,\ {\rm for}\ i=1,2,3.
\end{equation}

 It is
easy to see that $s_1s_2s_3\in \mathbb{Z}[t]^2,$ and that, for each $i$ and each prime $p\in \mathbb{Z}[t]$, we have
 \begin{equation}\label{ps123}
{\rm if}\ p|s_i\ {\rm then}\ p|s_js_k,\ {\rm where}\
i\neq j\neq k \neq i.
\end{equation}

Let $P\in E(\mathbb{Q}(t))\setminus\{O\}$. Then the first coordinate of $P$ is of the form
\begin{equation}\label{ikspe}
 x(P)=\frac{p(t)}{q(t)^{2}},\  {\rm with}\  p(t),q(t)\in \mathbb{Z}[t]\ {\rm coprime}
 \end{equation}
 (recall that $\mathbb{Z}[t]$ is an UFD). Therefore $p(t)-e_i(t)q^2(t)=s_i(P)\square_{\mathbb Z[t]},\ i=1,2,3$,
where  $\square_{\mathbb Z[t]}$ denotes a square of an element of
$\mathbb Z[t]$.
By this, \eqref{ps123} and the fact that $s_i$ are square-free, we deduce that
\begin{equation}\label{pse}
s_i|(e_j-e_i)(e_k-e_i),\mbox{ where }i\ne j\ne k\ne i
\end{equation}
for each $i$ (see also \cite{Huse}, Chapter 6, Proposition (4.1)). For example, a prime factor of $s_1$ is also a prime
factor of $s_2s_3$. Assume that it is a prime factor of $s_2$.
Then it is a prime factor of $(e_1-e_2)q^2(t)$, hence it is a
prime factor of $e_1-e_2$.

In Theorem \ref{main} we make a refinement of the method from
\cite{G-T}, Theorem 3.2. The proof is a modification of that
proof. Now we present the proof of Theorem \ref{main}.

\pff {\bf Theorem \ref{main}.}
 Note that the specialized curve is non-singular (see Lemma \ref{nonsingular} (i)). Let us prove that the specialization homomorphism is injective. Assume that the conditions of the theorem are satisfied and that $\sigma$ is not an injection. So there exists a point
$P\in E(\mathbb{Q}(t))\setminus \{O\}$ such that $\sigma(P)= O$. We will prove that it leads to a contradiction.
First we prove that $P\in 2E(\mathbb{Q}(t)).$  By \eqref{opetdupla}, it is equivalent to proving that
$s_i(t)=1$ for each $i=1,2,3$. Since $\sigma$ is injective on
the torsion part \cite[p. 272--273, proof of Theorem
III.11.4]{Sil3}, we may assume that $P\neq (e_i,0),\ i=1,2,3.$ By
$p(t)-e_k(t)q^2(t)=s_k(P)\square_{\mathbb Z[t]}$ and the fact that
$q(t_0)=0$, we get $p(t_0)=s_k(t_0)\square_{\mathbb Q}$. Since
$p(t_0)$ should be a non-zero rational square (recall that $q(t_0)=0$ and $p,q$ are coprime),
we see that $s_i(t_0)$ is a
rational square, for each $i=1,2,3$. We claim that $s_k(t)=1$ for
each $k=1,2,3$, i.e., that $P\in 2E(\mathbb{Q}(t))$.\\
Assume that $s_k(t)$ is non-constant for some $k$.
By the above discussion $s_k(t_0)$ is a rational square, which is in
contradiction with  condition (A) of the theorem (recall that by
\eqref{pse}, $s_k$ is a nonconstant square-free divisor of
$(e_i-e_k)\cdot (e_j-e_k)$ in  $\mathbb Z[t]$, with $i\neq j\neq k\neq i$).
 Therefore  $s_k(t)$ is constant for each $k$. Since $s_k(t)$ is
square-free in ${\mathbb Z[t]}$ and $s_k(t_0)$ is a rational
square, we see that $s_k(t)=1$, for each $k$. This proves that $P\in 2E(\mathbb{Q}(t))$.\\
We claim that there is $P_1\in E(\mathbb{Q}(t))$ such that $2P_1=P$ and $\sigma (P_1)=O$. Let
$P'_1\in E(\mathbb{Q}(t))$ be any point with $2P'_1=P$. Then
$2\sigma (P'_1)=O$, i.e., $\sigma(P'_1)$ is a $2$-torsion point on the specialized curve. Since $\sigma$ is injective on
the torsion points, there exists a $2$-torsion point $Q\in E(\mathbb{Q}(t))$  such that
$\sigma (Q)=\sigma (P'_1)$. Put $P_1=P'_1-Q$.  Then $2P_1=P$, especially $P_1\neq O$, and $\sigma(P_1)=O$. Note that $P_1$ is of infinite order.  Now
the procedure can be continued with $P_1$ instead of $P$, the contradiction. Therefore $P=O$, i.e., $\sigma$ is injective.
 \qed

\begin{remark}\label{slabiji} Condition (A) in Theorem \ref{main} is weaker then the
following condition\\
(A') For every nonconstant square-free divisor $h$ in  $\mathbb
Z[t]$ of
$$(e_1-e_2)\cdot (e_2-e_3)\cdot(e_3-e_1)$$
the rational number $h(t_0)$ is not a square in $\mathbb{Q}$.\\
For example, set $e_1=0,\ e_2=t,\ e_3=7t+1$. Then
$t_0:=\frac{1}{21}$ satisfies condition (A). Since
$\frac{1}{21}(6\cdot\frac{1}{21}+1)(7\cdot\frac{1}{21}+1)=(
\frac{2}{7})^2$, it does not satisfy condition (A').
\end{remark}

\begin{remark}\label{effective}
Let $\cal T$ denote the set of all integers
$t_0$ that satisfy Condition (A) from Theorem \ref{main}. Then there is an
effectively computable constant $c>0$ such that ${\cal T}\cap
[-c,c]\neq\emptyset$.
Namely, condition (A) in Theorem \ref{main} produces the equations  of the
form $z^2=h(t)$ for certain square-free polynomials $h$ over
$\mathbb{Z}$ of degree $d\geq 1$. If $d\leq 2$, the corresponding
curve has genus $0$, if $d=3$ or $4$ the genus is one, and if
$d\geq 5$ the curve is hyperelliptic with genus $\geq 2$. Recall
that  curves over $\mathbb{Q}$ of genus at least $1$ have only
finitely many integer points. Moreover, for elliptic and hyperelliptic curves,
there are explicit bounds for the height of integer points (\cite{Ba},  \cite{Bu}, Theorem 1; see also
 \cite{E-Si}, Theorem 1 b, for a bound of the number of integer points).
 If $d=1$ or $d=2$ then
the curve $z^2=h(t)$ may have finitely many or infinitely many integer
points. The case $d=1$ is straightforward, while the case $d=2$ reduces
to an estimating of the number of integer solutions for $Dz^2=t^2+B$ where $D$ is a square-free integer and $B$ a nonzero integer.
  The most demanding case is when $D\geq 2$. Then there is an
 effectively computable constant $c_1=c_1(D,B)$ such that $Dz^2=t^2+B$ has $\leq c_1 \tau (B)\log X$ integer solutions
 with $|t|,|z|\leq X$ for sufficiently large $X$, where $\tau(B)$ denotes the number of positive divisors of $B$ (see \cite{P-Z}, Lemma 3. for a more
precise estimation).
   \qed
\end{remark}

Now we sketch a proof of Theorem \ref{mainh}.  Assume  that $K$ is an arbitrary number field with  the ring of integers  $ \mathcal{O}_K$.
%\end{remark}
There exist at least one unique factorization domain $\mathcal R_K$, $\mathcal O_K\subset \mathcal R_K\subset K$ such that its group of units is finitely generated
(see for example \cite[p. 94, p. 127]{Knapp}).

\pff {\bf  Theorem \ref{mainh}.} Relations \eqref{esi}$-$\eqref{pse} remain valid after replacing
$\mathbb{Z}[t]$ by $\mathcal R_K[t]$. Now the proof is analogous to the proof of Theorem \ref{main}. Note that the theorem remains valid even if
we exclude the condition that $\mathcal R_K^{\times}$ is finitely generated.
However, this condition reduces the checking of Condition (C) from the Theorem
to checking of only a finitely many square-free divisors $h$ in
${\mathcal {R}}_K[t]$.
\qed
\medskip

It can be seen that there is a variant of Remark \ref{effective} for elliptic curves of the form \eqref{geneq}.
In the following remark we use another argument to prove that there are a
lot of rational integers $t_0$ satisfying condition (C) from Theorem
\ref{mainh}.
\begin{remark}
According to \cite{Sc}, Section 5, Definition 24,  Theorem 50 and Corollary 1,
for each $F\in\mathbb{C}[z,t]$ either:\\
(i) every congruence class
$\cal C$ in $\mathbb{Z}$ contains a congruence subclass ${\cal C}^{\ast}$ such that
for all $t_0\in{\cal C}^{\ast}$ the polynomial $F(z,t_0)$ has no zero in $K$, or\\
(ii)  $F$ viewed as a polynomial in $z$ has a zero in
$K(t)$.\\ By consecutively applying this to the polynomials
$F[z,t]:=z^2-h(t)$ above, we see that for each congruence class
$\cal C$ in $\mathbb{Z}$ there exists a congruence subclass $\cal C^{\ast}$
of $\cal C$, such that the conditions from Theorem \ref{mainh}  are
satisfied for all $t_0\in \cal C^{\ast}$.
\end{remark}

\section{Elliptic curves with exactly one rational $2$-torsion point}
\label{section3}

In this section we prove Theorem \ref{novi}, which  extends the criterion for injectivity from Theorem \ref{main} to elliptic curves having exactly one nontrivial rational $2$-torsion point. After that we sketch a proof of Theorem \ref{mainnovi}. Similarly as in Section \ref{section2} we first work over $\mathbb{Q}$ as the field of constants and after that we extend the consideration to arbitrary number fields.  Recall that
 $E:y^2=x^3+Ax^2+Bx+C,$  $A,B,C\in \mathbb{Z}[t]$ is a non-constant elliptic curve over $\mathbb{Q}(t)$, and that $D $ denotes the discriminant of the polynomial
$f(x):=x^3+Ax^2+Bx+C$. Since we assume that $E$ has exactly one $\mathbb{Q}(t)$-rational point we may work with the equation
\begin{equation}\label{tocnojedna}
E:y^2=x^3+Ax^2+Bx,\  A,B\in \mathbb{Z}[t],\ A^2-4B\notin \mathbb{Z}[t]^2,
\end{equation}
with $D=B^2(A^2-4B)$. Note that the factor $B$ corresponds to $e_1^2-(e+\bar e)e_1 +e\bar e$ from Theorem \ref{novi}, while the factor $A^2-4B$ corresponds to $(e-\bar e)^2$. The following condition concerning the discriminant $D$ and a rational number $t_0$ is a reformulation of the condition from Theorem \ref{novi}.\\

($\mathcal{A}$) For each factor $h$ of $B$ or of $A^2-4B$ in $\mathbb{Z}[t]$ if $h(t_0)$ is a square in $\mathbb{Q}$, then $h$ is  a square in $\mathbb{Z}[t]$.\\

\begin{lemma}\label{nonsingular}
Assume that  a rational number $t_0$ satisfies condition ($\mathcal{A}$).  Then:\\
(i) $D(t_0)\neq 0$.\\
(ii) The polynomial $f(x,t_0):= x^3+A(t_0)x^2+B(t_0)x$ has exactly one $\mathbb{Q}$-rational root.
\end{lemma}
\pf
  (i) Contrary, $D$ has a linear factor $h$ with $h(t_0)=0$, which contradicts condition ($\mathcal{A}$).\\
  (ii) Since $A^2-4B$ is not a square in $\mathbb{Z}[t]$, we conclude by ($\mathcal{A}$) that $A(t_0)^2-4B(t_0)$ is not a square in $\mathbb{Q}$.
  \qed
\par
 For the proof of Theorem \ref{novi} we will need some wellknown facts on $2$-isogeny (see, for example, \cite{Huse}, Chapter 4, Section 5). Let us define
 \begin{equation}\label{dualna}
\bar E:y^2=x^3-2Ax^2+(A^2-4B)x,
\end{equation}
and let $\bar D$ denote the discriminant of the polynomial $x^3-2Ax^2+(A^2-4B)x$. We have $\bar D=16B(A^2-4B)^2$. The map $\phi:E\rightarrow \bar E$ defined on points different from $O$ and $(0,0)$ by
$$\phi(P)=\left(\frac{y(P)^2}{x(P)^2},\frac{y(P)(x(P)^2-B)}{x(P)^2}\right),$$
gives rise to  an isogeny of degree two with kernel $\{O,(0.0)\}$ and with the dual isogeny defined by
$$\psi(\bar P)=\left(\frac{y(\bar P)^2}{4x(\bar P)^2},\frac{y(\bar P)(x(\bar P)^2-(A^2-4B))}{8x(\bar P)^2}\right).$$
  Further, $\bar P\in \bar E(\mathbb{Q}(t))$ different from $O$ and $(0,0)$ is in $\phi(E(\mathbb{Q}(t)))$ if and only if $x(\bar P)$ is a square from $\mathbb{Q}(t)$ (see \cite{Huse}, Chapter 4, Proposition 5.7  or \cite{SiTa}, p.83-85). From this it is easy to see that $P\in E(\mathbb{Q}(t))$ different from $O$ and $(0,0)$ is in $\psi(\bar E(\mathbb{Q}(t)))$ if and only if $x(P)$ is a square from $\mathbb{Q}(t)$.

\pff {\bf Theorem \ref{novi}.}
By Lemma \ref{nonsingular} (i), the specialized curve is non-singular, hence it is an elliptic curve over $\mathbb{Q}$. It remains to prove that the specialization homomorphism is injective. We may assume that $E$ is given by the equation of shape \eqref{tocnojedna}. Assume that $P\in E(\mathbb{Q}(t))$ is nontrivial with trivial specialization at $t_0$. We will see that this leads to the contradiction. Let us write $P=(x(P),y(P))$, where
\begin{equation}\label{skracen}
x(P)=\frac{p}{q^2}
\end{equation}
with coprime $p,q$ from $\mathbb{Z}[t]$, especially $q(t_0)=0$ and $p(t_0)$ is a square in $\mathbb{Q}$.\\ We claim that $P\in 2E(\mathbb{Q}(t))$. To prove this claim we first prove that
 $x(P)$ is a square in $\mathbb{Q}(t)$. It is enough to prove that $p$ is a square in $\mathbb{Z}[t]$. Contrary,  $p$ has a prime $\pi$ in $\mathbb{Z}[t]$ with an odd multiplicity. Since $p(p^2+Apq^2+Bq^4)$ is a square in $\mathbb{Z}[t]$, we see that $\pi$ is also a prime factor of $B$. Therefore, the square-free part $h$ of
$p$ in $\mathbb{Z}[t]$ has at least one prime factor and $h$ divides $B$. By condition ($\mathcal{A}$) of the Theorem $h(t_0)$ is not square in $\mathbb{Q}$. On the other side, since $p(t_0)$ is a square in $\mathbb{Q}$ and $q(t_0)=0$ we get that $h(t_0)$ is a square in $\mathbb{Q}$. This is the contradiction, hence $x(P)$ is a square in $\mathbb{Q}(t)$ (note that $-p$ is not a square in $\mathbb{Z}[t]$).
 From this we see that there exists a point $\bar P\in \bar E(\mathbb{Q}(t))$ such that $\psi (\bar P)=P$ (see the text after the proof of Lemma \ref{nonsingular}). It is easy to see that the specialization of $\bar E$  at $t_0$ is an elliptic curve and that $\bar P$  specializes to $O$ or to $(0,0)$. If $\bar P$ specializes to $(0,0)$ then $\bar P+(0,0)$ specializes to $O$, so we may assume that $\bar P$ specializes to $O$. Using the argument from the first part of the proof we conclude that $x(\bar P)$ is a square in $\mathbb{Q}(t)$, and further that there is a point $P_1\in E(\mathbb{Q}(t))$ such that $\bar P=\phi (P_1)$. Finally, we get $P=\psi (\bar P)=\psi(\phi(P_1))=2P_1$, as we claimed.
\par
 The rest of the  proof is analogous to the end of the proof of Theorem \ref{main}. We only have to note that, by Lemma \ref{nonsingular} (ii), the specialized curve $E_{t_0}$ has exactly one non-trivial $\mathbb{Q}$-rational $2$-torsion point.
\qed

We support the criterion from Theorem \ref{novi} by two examples.

\begin{example}\label{rang1}
Let $E$ be the elliptic curve over $\mathbb Q(t)$ given by the equation
$y^2=x^3+t^2x^2-x.$  Using the Tate-Shioda formula (see \cite{Sh2}, Corollary 5.3 and Lemma 10.1), one can find that it has rank $1$ over
 $\mathbb Q(t)$ with the point $P=(1,t)$ of infinite order. By  computing $mP$ for $m=2,3,...,12$
 we see that the specialization is injective for all rational $t_0$ except
$t=0,\pm 1$. This is not in collision with Theorem \ref{novi} because $t=0,\pm 1$ do not satisfy condition ($\cal A$).

\end{example}

\begin{example}\label{Bremner}
Let $E$ be  an elliptic curve $X$ over $\mathbb Q(t)$  given by the
equation
$y^2=x(x^2-2(5(2t^2-2t+1)(t^2-2t+2)-2(t^2-1)^2)x+25(2t^2-2t+1)^2(t^2-2t+2)^2)$
(see \cite{Br}, p. 551,  formula (14')).  By a detailed
$2$-descent analysis A. Bremner proved that the rank of $E$ over
$\mathbb Q(i)(t)$ is three, especially  the rank of $E$ over
$\mathbb Q(t)$ is as most three. He presented three independent
$\mathbb Q(t)$-rational points, which implies that the rank of $E$
over $\mathbb Q(t)$ is exactly three. To illustrate the criterion
from Theorem \ref{novi} let us note that the discriminant of $E$
equals to $$-2^8\cdot 5^4\cdot
(t-1)^2(t+1)^2(9t^4-30t^3+47t^2-30t+9)(t^2-2t+2)^4
(2t^2-2t+1)^4.$$ Further, $E$ has a nontrivial $2$-torsion point
$(0,0)$ and two $2$-torsion points conjugate over $\mathbb
Q(t,\sqrt{-9t^4+30t^3-47t^2+30t-9})$. Since $t_0=\frac{5}{2}$
satisfies condition ($\cal A$), and the specialized curve has rank
three over $\mathbb Q$, we conclude, by Theorem \ref{novi}, that
the rank of $E$ over $\mathbb Q(t)$ is at most three.
\end{example}

In the following remark we present two examples which show that
the criterion from Theorem \ref{novi} has not a direct extension to general
elliptic curves over $\mathbb Q(t)$ i.e. elliptic curves
$E:y^2=x^3+Ax^2+Bx+C,$ where $f(x):=x^3+Ax^2+Bx+C$ is irreducible
over $\mathbb Q(t)$. Let us consider the condition\\
($\cal A_1$) For each factor $h$ of $D$  in $\mathbb{Z}[t]$ if $h(t_0)$ is a square in $\mathbb{Q}$, then $h$ is  a square in $\mathbb{Z}[t]$.\\
It is clear that for polynomials $f$ having at least one nontrivial rational root, if $t_0$ satisfies condition ($A$) or condition ($\cal A$), then it satisfies condition ($\cal A_1$). For polynomials $f$ that are irreducible over $\mathbb Q(t)$ we consider the additional condition:
($\cal B$) The polynomial $f(x,t_0):=
x^3+A(t_0)x^2+B(t_0)x+C(t_0)$ is irreducible over $\mathbb Q$.

\begin{remark}\label{nijeopci}
\begin{itemize}

\item Let $E$ be the elliptic curve over $\mathbb Q(t)$ given by
the equation $y^2=x^3-x+t^2.$ It has rank $2$ over $\mathbb Q(t)$
with trivial torsion and generators any two of the points $(0,t),(1,t),(-1,t)$
(see \cite{Sh3}, Theorem $(A_2)$). Further its discriminant equals
to $16(4-27t^4)$ and the Galois group of the polynomial
$x^3-x+t^2$ is isomorphic to the symmetric group $S_3$. It is easy
to see that $t_0=\pm 1,\pm \frac 12$ satisfy conditions ($\cal A_1$) and ($\cal B$).
On the other side, one can check that the specializations at $\pm
1,\pm \frac 12$ are not injective.

 \item  Let $E$ be the elliptic
curve over $\mathbb Q(t)$ given by the equation $y^2=x^3-t^2x+1.$
It has rank $3$ over $\mathbb Q(t)$ with trivial torsion and
generators  $(0,1), (-1,t), (-t,1)$ (see  \cite{Tad},
Proposition 2.1 (iv)).  Its discriminant equals to $16(4t^6-27)$
and the Galois group of the polynomial $x^3-t^2x+1$ is isomorphic
to the symmetric group $S_3$. It is easy to see that $t_0=0$ satisfies condition ($\cal A_1$), while $t_0=\pm 1,\pm 2$
satisfy conditions ($\cal A_1$) and ($\cal B$). On the other side, one can check
that the specializations at $0,\pm 1\pm 2$ are not injective.

\end{itemize}

\end{remark}

Now we sketch a proof of Theorem \ref{mainnovi}. First note that the facts on $2$-isogeny that we used in the proof of Theorem \ref{novi} are valid for elliptic curves over arbitrary fields of characteristic zero.
Also condition ($\cal A$) can be extended directly to elliptic curves over $K(t)$ for a number field $K$ of class number $1$ (we only have to replace $\mathbb Z[t]$ by $\mathcal O_K[t]$). Now the proof of Theorem  \ref{novi} can be extended directly to this case (all we need is the unique factorization in $\mathcal O_K[t]$). It is clear that the later proof can be extended to arbitrary number fields $K$ (we only have to replace $\mathcal O_K[t]$ by $R_K[t]$).

\section{Application to an example by Mestre}
\label{section4}

In this section we apply our criterion to a family of quadratic twists of Mestre from the following example.

\begin{example}(\cite{Me}, \cite[Theorem 3.7]{R-S1}, \cite[Theorem 3]{Ste-T})\label{exam}
Let $a,b\in\mathbb Q$ such that $ab\ne 0$, let
$$g(t)=g^{a,b}(t)=-ab\cdot (t^2+1)\cdot (b^2(t^4+t^2+1)^3+a^3t^4(t^2+1)^2)$$
 and let $E=E^{a,b}$ be the elliptic curve over $\mathbb Q$ given by the equation
 \begin{equation}\label{eqjedn}
 y^2=x^3+ax+b.
 \end{equation}
 Then $E_g=E^{a,b}_g:y^2=x^3+ag(t)^2x+bg(t)^3$  has rank at least $2$ over $\mathbb Q(t)$, with two independent points $P=P_g^{a,b}$ and $Q=Q_g^{a,b}$
 with  coordinates
{\small{ \begin{equation}\label{P}
 P=\left(-\frac ba\frac{(t^2+t+1)(t^2-t+1)}{(t^2+1)}\cdot g(t),\frac 1{a^2(t^2+1)^2}\cdot g(t)^2\right)
 \end{equation}
 \begin{equation}\label{Q}
 Q=\left(-\frac ba\frac{(t^2+t+1)(t^2-t+1)}{t^2(t^2+1)}\cdot g(t),\frac 1{a^2t^3(t^2+1)^2}\cdot g(t)^2\right).
 \end{equation}
  }}
\end{example}
\medskip

Similarly as in \cite{Ste-T}, Section 4, let $C$ be a smooth
complete model over $\mathbb{Q}$ of the curve given by $s^2=g(t)$.
Then for each point $T=(x(t),y(t))\in E_g(\mathbb{Q}(t))$ there is
a morphism $\phi_T:C\rightarrow E$ defined by
$\phi_R(t,s)=(\frac{x(t)}{s^2},\frac{y(t)}{s^3})$. This gives rise
to an isomorphism between ${\rm Mor}_{\mathbb{Q}}(C,E)$ modulo
constant morphisms and $E_g(\mathbb{Q}(t))$ modulo torsion. The
degree $\deg \phi_T$ equals to $\deg \frac{x(t)}{g(t)}$. Since
$\deg (2\phi)=4\deg \phi$ for each $\phi\in {\rm Mor}(C,E)$ the
mapping $T\mapsto \frac{1}{2}\deg \phi_T$ is the canonical
 height on $E_g(\mathbb{Q}(t))$ (for a direct proof see \cite{G-L}).
Let $\langle\ ,\ \rangle$ denote the corresponding canonical
bilinear form on $E_g(\mathbb{Q}(t))\times E_g(\mathbb{Q}(t))$,
i.e.
$$\langle T,S\rangle=\frac{1}{2}(\deg \phi_{T+S}-\deg \phi_T-\deg \phi_S),\ {\rm for\ each}\ T,S\in
E_g(\mathbb{Q}(t)),$$ especially $\deg \phi_T=\langle T,T\rangle,\
{\rm for\ each}\ T.$
This approach works for general nonconstant twists (quadratic or else) of constant elliptic curves. See \cite[\S 4.2, \S 4.4]{Pa} for some further aspects of this approach which is specific for the Mestre family.\\
Since $\deg \phi_P=\deg \phi_Q=4$ and $\deg \phi_T\geq 4$ for all nontorsion points $T$, it is intuitively clear that
$P,Q$ generate a maximal nontorsion subgroup of rank two in $E_g(\mathbb{Q}(t))$. In the following lemma we give a precise formulation and a proof.
\begin{lemma}\label{GLT}
Let $E_g=E_g^{a,b}$, $P=P^{a,b}_g, Q=Q^{a,b}_g$ be as in  Example \ref{exam}.
Then:\\
(i) $P,Q$ are free generators of each subgroup of rank two in
$E_g(\mathbb{Q}(t))$ containing $P,Q$. \\
(ii) Assume that there exists a $t_0\in \mathbb{Q}$ for which the corresponding specialization
homomorphism $\sigma_{t_0}$ is injective and the rank of the
specialized elliptic curve over $\mathbb Q$ is two. Then the rank of
$E_g(\mathbb{Q}(t))$ is two and $P,Q$ are its free generators.
\end{lemma}
\pf 
(i) Let $M$ be a torsion-free subgroup of $E_g(\mathbb{Q}(t))$ of rank
two containing $P,Q$ and let $T\in M$ be a nontorsion point.
 Then there is a nontrivial relation
\begin{equation*}
kT=mP+nQ,\ m,n,k\in \mathbb{Z}.
\end{equation*}
We may assume that  $k>0$. By consecutive adding or subtracting $kP$ or $kQ$, it
leads to
\begin{equation*}
kT_1=m'P+n'Q,\ {\rm with}\ -\frac{k}{2}\leq m',n'\leq \frac{k}{2}.
\end{equation*}

From this we get
\begin{equation*}
k^2\deg \phi_{T_1}=m'^2\deg \phi_P+n'^2\deg \phi_Q+2m'n'\langle P,Q\rangle,
\end{equation*}
which provides an upper bound for $\deg \phi_{T_1}$:
\begin{equation}\label{visina}
\deg \phi_{T_1}\leq \frac{\deg \phi_P+\deg \phi_Q+2|\langle P,Q\rangle|}{4}.
\end{equation}
Since $\deg \phi_P=\deg \phi_Q=4$, $\deg \phi_{P+Q}\leq 8$ and $\deg \phi_{P-Q}\leq 8$ for each $a,b$, we conclude, by the parallelogram law that $\deg \phi_{P+Q}=\deg \phi_{P-Q}=8$,
hence $ \langle P,Q\rangle=0$. By \eqref{visina} we get $\deg \phi_{T_1}\leq 2$. We claim that $\deg \phi_{T_1}=1$ or $\deg \phi_{T_1}=2$ is impossible.  Contrary,
 $x=x(T_1)/g(t)=\frac{\alpha(t)}{\beta(t)}$, where $\alpha,\beta$ are nonzero polynomials over $\mathbb{Q}$ of degree at most $2$ and at least one of them is non-constant. Plugging in $g(t)y^2=x^3+ax+b$, we get that there is a nonzero polynomial $w$ over $\mathbb{Q}$ of degree at most $6$, such that $w(t)\beta(t)g(t)$ is a square in $\mathbb{Q}[t]$. It is impossible since $g(t)$ is squarefree with degree $14$. Therefore $\deg \phi_{T_1}=0$, hence $T_1$ is torsion, i.e. $T$ is a $\mathbb{Z}$-linear combination of $P,Q$ and torsion points.\\
 (ii) Directly from (i).
\qed

 Let $K$ denote the splitting  field
 of the cubic polynomial $f(x)=x^3+ax+b$,
thus $K$ is Galois. It is well known that either $K=\mathbb{Q}$, $ K$ is a quadratic field over $\mathbb{Q}$,  $K$ is a cubic field over $\mathbb{Q}$ with cyclic Galois group, or $K$ is a sextic field over $\mathbb{Q}$ with the Galois group isomorphic to the symmetric group $\mathbf{S}_3$.

In the sequel we illustrate Theorem \ref{mainh} on two concrete examples from Mestre's family to get the rank and prove that given points are  free generators, we describe the corresponding algorithm. In fact, we performed a more extensive calculation for various number fields $K$ of class number one (including $\mathbb Q$) with various Galois groups, and for several elliptic curves for a number field $K$ of class number two. In all cases we get that the rank is two and proved that the given points are  free generators. One can find these calculations in \cite[Section 5]{GTarxiv}. Two examples are presented here in details, one when the splitting field $K$  is a sextic field with class number one (see Example \ref{corall}), the other  when $K$ has class number two (see Example \ref {corh2}).
Calculations are performed using a variety of packages: GP/Pari \cite{Pari}, MAGMA
\cite{MAGMA}, mwrank \cite{mwrank}.
Let us sketch the algorithm. After fixing a concrete value of $(a,b)$ we chose a rational number $t_0$ such that:
\begin{itemize}
\item $\sigma_{t_0}$ is an  injection (by using  Theorem
\ref{mainh}).
\item   $t_0$ is such that the specialized curve $E^{a,b}_g(t_0)$ is of rank 2 over $\mathbb Q$
(which is calculated with {\tt mwrank} and Magma's command  {\tt MordellWeilShaInformation}). To avoid extensive calculations, before calculating the rank we checked if the root number was one.
\end{itemize}
Now  Lemma \ref{GLT}(ii)  is applied to conclude that $E_g=E^{a,b}_g$
over $\mathbb Q(t)$ has actually free generators the two points
$P,Q$ from Example \ref{exam}.
If $K$ properly contains the field of rational numbers Theorem \ref{mainh} gives the injectivity of $\sigma_{t_0}$ as a homomorphism from  $E_g(K(t))$ to $E_g(t_0)(K)$, so we have to  look at the restriction of $\sigma_{t_0}$ to $E^{a,b}_g(\mathbb
Q(t))$, which is injective, too.

We first present an example in which the splitting  field $K$
 of the cubic polynomial $f(x)=x^3+ax+b$  has class number one with the Galois group the symmetric group  $\mathbf{S}_3$.
\begin{example} \label{corall}
Let $E=E^{1,1}$ be the elliptic curve over $\mathbb Q$ given by equation \eqref{eqjedn} with $a=b=1$
$$E^{1,1}:y^2=x^3+x+1.$$
Then the elliptic curve $E_g$   has rank
two over $\mathbb Q(t)$, with free generators the two points $P=P^{1,1}$ and
$Q=Q^{1,1}$ listed in  Example \ref{exam}.
 The splitting field $K$  of the polynomial $x^3+x+1$ is the sextic field $\mathbb Q(q)$ of class number one, generated by  the algebraic number $q$ defined as a root of the polynomial $x^6 + 78x^4 + 324x^3 + 1521x^2 + 12636x + 64219$.  The two  fundamental units of $K$ are
$$\frac 4{245805}q^5 - \frac{169}{737415}q^4 + \frac{52}{49161}q^3 - \frac{7097}{737415}q^2 - \frac{8728}{49161}q - \frac{13156}{737415}, $$ $$\frac 2{49161}q^5 - \frac{169}{294966}q^4 + \frac{130}{49161}q^3 - \frac{7097}{294966}q^2 + \frac{5521}{98322}q - \frac{6578}{147483}.$$
 Further,
{\scriptsize{
$$e_1(t)=-\left(-\frac 2{35115}q^5 + \frac{169}{210690}q^4 - \frac{26}{7023}q^3 + \frac{7097}{210690}q^2 + \frac{1705}{14046}q + \frac{6578}{105345}\right)\cdot g^{1,1}(t),$$
$$e_2(t)=-\left(\frac 4{245805}q^5 - \frac{169}{737415}q^4 + \frac{52}{49161}q^3 - \frac{7097}{737415}q^2 - \frac{8728}{49161}q - \frac{13156}{737415}\right)\cdot g^{1,1}(t),$$
$$e_3(t)=-\left(\frac 2{49161}q^5 - \frac{169}{294966}q^4 + \frac{130}{49161}q^3 - \frac{7097}{294966}q^2 + \frac{5521}{98322}q - \frac{6578}{147483}\right)\cdot g^{1,1}(t).$$
}}
If we choose $t_0=3$, it is easy to see  that if $h(t)$ is a nonconstant divisor of one of the
rad$_{\mathcal O_K[t]}((e_k(t)-e_i(t))\cdot(e_k(t)-e_j(t))$ ($i,j,k\in\{1,2,3\}$ different)  in $\mathcal O_K[t]$,
then $h(3)$ is not a square in $ K.$ Thus Theorem \ref{mainh} is
satisfied for $K=\mathbb Q(q)$. We conclude that the
specialization homomorphism $\sigma_{3}:E^{1,1}_g(K(t))\rightarrow
E^{1,1}_g(3)(K)$ is an injection, hence its restriction to ${E_g(\mathbb
Q(t))}$
is also an injection.
The rank of
$E^{1,1}_g(3)(\mathbb Q)$ is $2$. Thus we have shown  by using Lemma \ref{GLT}(ii)
that $E_g=E^{1,1}_g$   has rank
two over $\mathbb Q(t)$, with free generators $P,Q$.
\qed
\end{example}
In the  following example the splitting field of $f(x)=x^3+ax+b$
is  a field of class number two, specifically $K=\mathbb
{Q}(\sqrt{-5})$.  We apply Theorem \ref{mainh} and the result
comment from the viewpoint of Theorem \ref{novi}. First we have to
choose an adequate $\mathcal R_K=\mathcal R_{\mathbb
Q(\sqrt{-5})}$ as in Theorem \ref{mainh}. The choice of $\mathcal
R_K$ in general isn't unique. To reduce the number of divisors in
the condition of Theorem \ref{mainh} to a finite number, we choose
$\mathcal R_K$ with finitely generated group of units.

By the example in \cite[p.129]{Knapp} we know that we can choose for $\mathcal R_{\mathbb Q(\sqrt{-5})}$  the principal ideal domain (PID) which is the localization of $\mathcal O_K=\mathcal O_{\mathbb Q(\sqrt{-5})}$ by the multiplicative set $S=\{1,2,2^2,2^3,2^4,2^5,\ldots\}$, where the group of units is generated by $\mathcal O_{\mathbb Q(\sqrt{-5})}^{\times}=\{\pm 1\}$ and $2$. So $K=\mathbb Q(\sqrt{-5})$, $O_K=\mathbb Z[\sqrt{-5}]$ and $\mathcal R_K=S^{-1}\mathcal O_K$. For each ideal $I$ of $\mathcal O_K$ let us define $I_S:=S^{-1}I$. Then $I_S$ is an ideal of $\mathcal R_K$ and $I_S$ is proper if and only if $I\cap S=\emptyset$. The non-zero prime ideals of $\mathcal R_K$ are exactly $S^{-1}I$ where
$I$ goes through non-zero prime ideals of $\mathcal O_K$ different from $\mathcal{P}=(2,1-\sqrt{-5})$.
Since $K$ is the quotient field of the unique factorization domain $\mathcal R_K$, thus we can obtain the irreducible nonconstant factors of a polynomial in $\mathcal R_K[t]$ by observing the factorization in $K[t]$.

\begin{example} \label{corh2}
Let $E=E^{2,12}$ be the elliptic curve over $\mathbb Q$ given by equation \eqref{eqjedn} with $(a,b)=(2,12)$. Then the elliptic curve $E_g=E^{2,12}_g$
 has rank
two over $\mathbb Q(t)$, with free generators the two points $P=P^{2,12}$ and
$Q=Q^{2,12}$ listed in  Example \ref{exam}.\\
The splitting field of the polynomial $x^3+2x+12$ is $K=\mathbb Q(\sqrt{-5})$ which is of class number two. In this case we have
$$g(t)=g^{2,12}(t)=-2^6\cdot 3\cdot  (t^2+1)(3t^4+2t^2+2)(3t^4+4t^2+3)(2t^4+2t^2+3).$$
Thus  we look at the elliptic curve $E_g=E^{2,12}_g$ over $K(t)$.
 Further from the coefficients we get that the discriminant is equal to
\noindent
$(e_1(t)-e_2(t))^2\cdot(e_1(t)-e_3(t))^2\cdot (e_2(t)-e_3(t))^2=$\\
$$-2^{40}\cdot 3^6\cdot 5\cdot 7^2\cdot (t^2+1)^6(3t^4+2t^2+2)^6(3t^4+4t^2+3)^6(2t^4+2t^2+3)^6,$$
which we have to factor into irreducibles in $\mathcal R_K[t]$, only the radical is of importance to get the square-free factors in $\mathcal R_K[t]$.
One shows $\sqrt{-5},$ $1\pm \sqrt{-5},$  $3\pm \sqrt{-5}$ are irreducible elements  in the principal ideal domain $\mathcal{R}_{K}=\mathcal{R}_{\mathbb Q(\sqrt{-5})}$ described above. $2$ and $-1$ are invertible elements in $\mathcal R_K$ and $\mathcal R_K[t]$. We also have $3=\frac 12(1+\sqrt{-5})(1-\sqrt{-5}),\ 2\pm \sqrt{-5}=-\frac 12(1\mp\sqrt{-5})^2.$

We first factor in $K[t]$ to obtain the factorization in $\mathcal R_K[t]$.
Now it is easy to see that the radical in $\mathcal R_K[t]$ of $(e_1(t)-e_2(t))\cdot(e_1(t)-e_3(t))\cdot (e_2(t)-e_3(t))$ factors into irreducible elements in the UFD $\mathcal R_K[t]$ as
$\mbox{rad}_{\mathcal R_K[t]}[(e_1(t)-e_2(t))\cdot(e_1(t)-e_3(t))\cdot (e_2(t)-e_3(t))]=$

$$= \sqrt{-5}\cdot (1+\sqrt{-5})\cdot (1-\sqrt{-5})\cdot (3+\sqrt{-5})\cdot (3-\sqrt{-5})\cdot (t^2 + 1)\cdot$$
$$\cdot   \left(t^2 + \frac{1+\sqrt{-5}}2\right)  \left(t^2 +\frac {  1-\sqrt{-5}}2\right)    \left(\frac  {1-\sqrt{-5} }2t^2 +  1\right)  \left(\frac  {1+\sqrt{-5} }2t^2 +1\right)\cdot$$ $$\cdot \left((1+\sqrt{-5})t^2 - (1-\sqrt{-5})\right)   \left((1-\sqrt{-5} )t^2 -(1+\sqrt{-5})\right) .$$
So we obtain all nonconstant square-free divisors of $(e_1(t)-e_2(t))\cdot(e_1(t)-e_3(t))\cdot (e_2(t)-e_3(t))$ in $\mathcal R_K[t]$,
we had to take into account  $-1$ and $2$, the two generators of
the group of units in $\mathcal R_K$. If we choose $t_0=4$ using
analogous arguments as in  the previous example we conclude
$E^{2,12}_g$ over $\mathbb Q(t)$ has rank two and free generators
$P,Q$.\\
Note that $t_0=4$ satisfies condition ($\cal A$) of Theorem
\ref{novi}, which confirms that the specialization at $t_0=4$ is
injective.

\end{example}

\medskip

{\small{
{\footnotesize \noindent FACULTY OF CHEMICAL ENGINEERING AND TECHNOLOGY, UNIVERSITY OF ZAGREB, Maruli\'cev trg 19, 10000 Zagreb, Croatia \\
{\em E-mail address}, I. Gusi\' c: {\tt igusic@fkit.hr}}

{\footnotesize \noindent DEPARTMENT OF ANALYSIS AND COMP. NUMBER THEORY, GRAZ UNIVERSITY OF TECHNOLOGY, Steyrergasse 30, 8010 Graz, Austria\\
 {\em  E-mail address}, P. Tadi\' c: {\tt petra.tadic.zg@gmail.com
 }
}}

\end{document}